\documentclass[11pt,reqno,a4paper]{amsart}
 \usepackage{amsgen, amstext,amsbsy,amsopn, amsthm, amsfonts,amssymb,amscd,amsmath,euscript,enumerate,url,verbatim,calc,tikz}
\usepackage{hyperref}
\usepackage{MnSymbol}
\usepackage{pgfkeys}
\usepackage{mathtools}
\usepackage{etex}
 \RequirePackage{etex}
\usepackage[left=1.2in,right=1.2in,bottom=1.5in]{geometry}

\usepackage{thmtools, thm-restate}
\usepackage{pst-node}

\usepackage{mathrsfs}
\usepackage{tikz-cd}

\usepackage{tikz}
\usepackage{tikz-network}
\usepackage{eqnarray,amsmath}

\let\emptyset\varnothing

\def\multiset#1#2{\ensuremath{\left(\kern-.3em\left(\genfrac{}{}{0pt}{}{#1}{#2}\right)\kern-.3em\right)}}

\usetikzlibrary{arrows}

 \usepackage{latexsym}
 \usepackage{graphics}
 \usepackage{color}
\usepackage{lastpage}
\usepackage{fancyhdr}
\usepackage{multirow}
\usepackage[T1]{fontenc}
\usepackage[utf8]{inputenc}
\allowdisplaybreaks
\usepackage{graphicx}
\graphicspath{ {F:/IMAGES/} }

\makeatletter
\def\oversortoftilde#1{\mathop{\vbox{\m@th\ialign{##\crcr\noalign{\kern3\p@}%
      \sortoftildefill\crcr\noalign{\kern3\p@\nointerlineskip}%
      $\hfil\displaystyle{#1}\hfil$\crcr}}}\limits}

\def\sortoftildefill{$\m@th \setbox\z@\hbox{$\braceld$}%
  \braceld\leaders\vrule \@height\ht\z@ \@depth\z@\hfill\braceru$}

\makeatother

  \newcommand{\Ass}{\operatorname{Ass}}

\newcommand{\proset}{\,\mathrel{\lower 4pt\hbox{$\scriptscriptstyle/$}
\mkern -14mu\subseteq }\,} 

 \newtheorem{theorem}{Theorem}[section]
 \newtheorem{corollary}[theorem]{Corollary}
 \newtheorem{lemma}[theorem]{Lemma}
 \newtheorem{proposition}[theorem]{Proposition}

\usepackage{amsmath}

\usepackage{hyperref}

 \theoremstyle{definition}
 
 \newtheorem{remark}[theorem]{Remark}
 \newtheorem{definition}[theorem]{Definition}

 \newtheorem{example}[theorem]{Example}

\title[v-Number and the Newton-Okounkov region]{Asymptotic v-number of graded families of ideals  and the Newton-Okounkov region}
\author{Mousumi Mandal and Partha Phukan}

\thanks{AMS Classification 2020: 13A02, 13A15, 13F20, 13F55, 05E40}
\thanks{Key words and phrases: v-number, regularity, Newton-Okounkov region, multiplicity, integral closure}
\address{Department of Mathematics, Indian Institute of Technology Kharagpur, 721302, India}\email{mousumi@maths.iitkgp.ac.in}
\address{Department of Mathematics, Indian Institute of Technology Kharagpur, 721302, India}\email{p.partha.24@kgpian.iitkgp.ac.in}

\begin{document}

\maketitle

\begin{abstract}
In this paper, we prove that for Noetherian graded families $\mathcal{I} = \{I_k\}_{k \ge 0}$ of homogeneous ideals in a Noetherian $\mathbb{N}$-graded Noetherian domain, $\lim\limits_{k \to \infty} \frac{\mathrm{v}(I_k)}{k}$ exists, 
and is given by $\frac{\alpha(I_r)}{r}$ for some $r \ge 1$, where $\alpha(I)$ denotes the initial degree. Extending these results to integral closures, we show that
$
\lim\limits_{k\to\infty}\frac{\mathrm{v}(\overline{I_k})}{k} =  \lim\limits_{k\to\infty}\frac{\alpha(\overline{I_k})}{k}=\lim\limits_{k\to\infty}\frac{\mathrm{v}(I_k)}{k}=\lim\limits_{k\to\infty}\frac{\alpha(I_k)}{k}  
$. 
For a graded family of monomial ideals in a polynomial ring, we provide a combinatorial interpretation of these limits via Newton--Okounkov regions $\Delta(\mathcal{I})$.
This connection is further generalized to arbitrary homogeneous ideals using good valuations. We also establish that both $\operatorname{reg}(I_k)$ and $\mathrm{v}(I_k)$ are eventually quasi-linear functions of $k$ for any Noetherian graded family. 
For a stable monomial ideal $I$ we show that $\mathrm{v}(I) < \operatorname{reg}(I)$. Finally, for a zero-dimensional homogeneous ideal $I$ in a polynomial ring $S$, we prove that $\mathrm{v}(I) < e(S/I)$, where $e(S/I)$ denotes the  multiplicity.
\end{abstract}

\section{Introduction}

Let $R = \bigoplus_{d \geq 0} R_d$ be an $\mathbb{N}$-graded Noetherian ring with $R_0$ be an infinite field, and let $I \subseteq R$ be a proper homogeneous ideal. For any associated prime $\mathfrak{p} \in \operatorname{Ass}(I)$, there exists a homogeneous element $f \in R$ such that $(I:_R f) = \mathfrak{p}$. This elementary observation leads to the following invariant.

\begin{definition}
For a proper homogeneous ideal $I \subseteq R$, the {\it $\mathrm{v}$-number} of $I$, denoted $\mathrm{v}(I)$, is defined as
\[
\mathrm{v}(I) := \min\{ d \geq 0 \mid \exists \, f \in R_d \text{ and } \mathfrak{p} \in \operatorname{Ass}(I) \text{ such that } (I:_R f) = \mathfrak{p} \}.
\]
For each associated prime $\mathfrak{p} \in \operatorname{Ass}(I)$, one may also consider the {\it local $\mathrm{v}$-number} of $I$ at $\mathfrak{p}$, denoted by $\mathrm{v}_{\mathfrak{p}}(I)$, defined as
\[
\mathrm{v}_{\mathfrak{p}}(I) := \min\{ d \geq 0 \mid \exists \, f \in R_d \text{ such that } (I:_R f) = \mathfrak{p} \}.
\]
\end{definition}

The concept of the \(\mathrm{v}\)-number was first introduced by Cooper et al. in \cite{Cooper20} to study the asymptotic behaviour of the minimum distance function of projective Reed–Muller type codes. In combinatorics and graph theory, Jaramillo and Villarreal \cite{JaramilloVillarreal21} established a connection between the \(\mathrm{v}\)-number of edge ideals of clutters and the independent domination number of the clutter itself, thereby providing a combinatorial interpretation for the \(\mathrm{v}\)-number of square-free monomial ideals.

This paper focuses on studying the asymptotic behavior of the $\mathrm{v}$-number for Noetherian graded families of homogeneous ideals in an $\mathbb{N}$-graded Noetherian domain $R$. We aim to establish analogues for the $\mathrm{v}$-number of results from the paper ~\cite{HaNguyenNguyen25} by H. T. H\`a et al. In \cite[Theorem~3.4]{HaNguyenNguyen25}, they proved that if $\mathcal{I}=\{I_k\}_{k\geq 0}$ is a Noetherian graded family of homogeneous ideals in a reduced standard graded $K$-algebra $R$ with $\operatorname{Ann}_R(I_k)=(0)$ for all $k\geq 1$, then
\[
\lim_{k\to\infty}\frac{\operatorname{reg}(I_k)}{k}
=
\lim_{k\to\infty}\frac{d(I_k)}{k}
=
\lim_{k\to\infty}\frac{\operatorname{reg}(\overline{I_k})}{k}
=
\lim_{k\to\infty}\frac{d(\overline{I_k})}{k},
\]
where $d(J)$ denotes the maximal degree of a minimal homogeneous generating set of $J$. Motivated by this result, we ask whether an analogous statement holds for the $\mathrm{v}$-number, with the maximal generating degree replaced by the initial degree. We answer this question affirmatively for Noetherian graded families of homogeneous ideals in an $\mathbb{N}$-graded Noetherian domain $R$ (Theorem~\ref{mainresultofhomogeneousideals}). In \cite{FicarraSgroi23}, Ficarra and Sgroi proved that if \(S\) is a standard graded polynomial ring over a field \(K\) and \(I\) is a homogeneous ideal of \(S\), then
\(
\lim\limits_{k \to \infty} \frac{\mathrm{v}(I^k)}{k} = \alpha(I),
\) where \(
\alpha(I)
=
\min\{\deg(f)\mid 0\neq f\in I , \text{ is homogeneous}\}
\). Further generalizations and refinements of these results were considered in \cite{VanmathiSarkar25},\cite{FiorindoGhosh24},\cite{FiorindoGhosh},\cite{GhoshNanduriPramanik26},\cite{GhoshPramanik25},\cite{{KotalSaha24}}. In \cite{KumarNanduriSaha25}, Kumar et al. proved that \(\lim\limits_{k \to \infty} \frac{\mathrm{v}(I_k)}{k}\) exists for a Noetherian graded filtration of homogeneous ideals in \(R\). In this paper, we generalize this result to a Noetherian graded family of homogeneous ideals.

\begin{theorem}[Theorem {\ref{Proposition 1.4}}]
Let \(\mathcal{I} = \{I_k\}_{k \geq 0}\) be a Noetherian graded family of homogeneous ideals in a Noetherian \(\mathbb{N}\)-graded domain \(R\). Then
\[
\lim_{k \to \infty} \frac{\operatorname{v}(I_k)}{k}
\]
exists, and there exists \(r \geq 1\) such that
\[
\lim_{k \to \infty} \frac{\operatorname{v}(I_k)}{k} = \frac{\alpha(I_r)}{r}.\]
\end{theorem}

We now turn our attention to the integral closure of a Noetherian graded family of homogeneous ideals in \(R\), where we proved our main theorem.

\begin{theorem}[Theorem {\ref{mainresultofhomogeneousideals}}]
Let $\mathcal{I}=\{I_k\}_{k \ge 0}$ be a Noetherian graded family of homogeneous ideals in a Noetherian  $\mathbb{N}$-graded domain $R$. Then
\[
\lim_{k\to\infty}\frac{\mathrm{v}(\overline{I_k})}{k}
   =\lim_{k\to\infty}\frac{\alpha(\overline{I_k})}{k}
   =\lim_{k\to\infty}\frac{\mathrm{v}(I_k)}{k}
   =\lim_{k\to\infty}\frac{\alpha(I_k)}{k}.
\]
\end{theorem}

We now aim to interpret the asymptotic \(\mathrm{v}\)-number of a Noetherian graded family \(\{I_k\}_{k \ge 0}\) from a combinatorial perspective. To this end, we employ the framework of Newton--Okounkov regions $\Delta(\mathcal{I})$ associated with the family. This construction extends the classical notion of Newton--Okounkov bodies, which originated in the work of Okounkov \cite{Okounkov96,Okounkov03} and was subsequently developed by Lazarsfeld and Musta\c{t}\u{a}~\cite{LazarsfeldMustata09}, Kaveh and Khovanskii~\cite{{KavehKhovanskii12},{KavehKhovanskii14}}, and Cutkosky~\cite{ Cutkosky13, Cutkosky14}. The theory has attracted considerable attention and found numerous applications across diverse areas of mathematics.

\begin{theorem}[Theorem {\ref{analogue}}]

Let $\mathcal{I}=\{I_k\}_{k \ge 0}$ be a Noetherian graded family of monomial ideals in $S = K[x_1, \dots, x_n]$. Then
\[
\lim_{k\to\infty}\frac{\mathrm{v}(\overline{I_k})}{k}
   =\lim_{k\to\infty}\frac{\alpha(\overline{I_k})}{k}
   =\lim_{k\to\infty}\frac{\mathrm{v}(I_k)}{k}
   =\lim_{k\to\infty}\frac{\alpha(I_k)}{k}
   =\lambda(\Delta(\mathcal{I})),
\]
where $\lambda(\Delta(\mathcal{I})):=\min\{|\mathbf{v}|:\mathbf{v}\text{ is a vertex of }\Delta(\mathcal{I})\}$.
\end{theorem}

Next, we generalize this result from monomial ideals to homogeneous ideals using the notion of a good valuation.

\begin{theorem}[Theorem \ref{NORforvaluations}]

Let $\mathcal{I} = \{I_k\}_{k \ge 0}$ be a Noetherian graded family of homogeneous ideals, 
and let $v$ be a good valuation that respects the monomials in $S = K[x_1, \dots, x_n]$. 
Let $\Delta(\mathcal{I})$ be the Newton--Okounkov region of 
$\mathcal{I}$ defined by $v$.  Suppose that the limiting region $\mathcal{C}(\mathcal{I})$ is a polyhedron. Then
\[
\lim_{k \to \infty} \frac{\mathrm{v}(\overline{I_k})}{k} 
= \lim_{k \to \infty} \frac{\alpha(\overline{I_k})}{k}
= \lim_{k \to \infty} \frac{\mathrm{v}(I_k)}{k}
= \lim_{k \to \infty} \frac{\alpha(I_k)}{k}
= \lambda(\Delta(\mathcal{I})),
\]
where $\lambda(\Delta(\mathcal{I})) := \min\{\,|\mathbf{v}| : \mathbf{v} \text{ is a vertex of } \Delta(\mathcal{I})\,\}$.

\end{theorem}

A classical result in commutative algebra, proven independently by Cutkosky et al. \cite{CutkoskyHerzogTrung99} and Kodiyalam \cite{Kodiyalam00}, states that for a homogeneous ideal \(I\) in a polynomial ring over a field, the Castelnuovo--Mumford regularity \(\operatorname{reg}(I^n)\) is eventually a linear function of \(n\). Recent attention has turned to the asymptotic behaviour of the Vasconcelos invariant. Conca \cite{Conca24}, and independently Ficarra--Sgroi \cite{FicarraSgroi23}, proved that \(\operatorname{v}(I^n)\) is also eventually linear in \(n\), with the leading coefficient equal to  \(\alpha(I)\). In a subsequent work, Ficarra and Sgroi \cite{FicarraSgroi24} extended this result to a much broader setting. They considered a Noetherian graded filtration \(\mathcal{I} = \{I_{k}\}_{k\ge 0}\) of a Noetherian \(\mathbb{N}\)-graded domain \(R\) and proved that  \(\mathrm{v}(I_{k})\) is an eventually quasi-linear function in \(k\). This generalization unifies the previous results and applies to many important families of ideals, including ordinary powers, symbolic powers, and integral closures of powers. We now generalize these results to the setting of Noetherian graded families. The following theorems establish the quasi-linearity of the regularity and v-number for such families.

\begin{theorem}[Theorem \ref{reg/vnumberisquasilinear}]
    Let $\mathcal{I} = \{ I_k \}_{k \geq 0}$ be a Noetherian graded family of homogeneous ideals in a standard $\mathbb{N}$-graded Noetherian domain $R$. Then $\operatorname{reg}(I_k)$ is eventually a quasi-linear function in $k$.
        
\end{theorem}

\begin{theorem}[Theorem \ref{vnumberisquasilinearforhomogeneous}]
Let \(\mathcal{I} = \{I_k\}_{k \geq 0}\) be a Noetherian graded family of homogeneous ideals in a  $\mathbb{N}$-graded Noetherian domain \(R\). Then \(\mathrm{v}(I_k)\) is eventually a quasi-linear function in \(k\). 

\end{theorem}

The \(\mathrm{v}\)-number has been examined from several angles in the literature. One key observation is that it gives a lower bound for the Castelnuovo–Mumford regularity for various classes of ideals \cite{VanmathiSarkar25, AmbhoreSahaSengupta24,BiswasMandal24,Ficarra26,FicarraMarques25,JaramilloVillarreal21, KumarNanduriSaha25, Saha24, SahaSengupta22}. The comparison between the \(\mathrm{v}\)-number and regularity has drawn the attention of many researchers. In this paper, we prove that for stable monomial ideals, the inequality \(\mathrm{v}(I) < \operatorname{reg}(I)\) holds strictly.

\begin{corollary}[Corollary \ref{v(I)<reg(I)forstablemonomialideals}]

Let $I$ be a stable monomial ideal (strongly stable, lexsegment, universal lexsegment, respectively) in $S = K[x_1, \dots, x_n]$. Then \[\mathrm{v}(I) < \operatorname{reg}(I).\]

\end{corollary}

\noindent While much of the existing literature has focused on the relationship between the v-number and the Castelnuovo–Mumford regularity, in this paper we explore a different but equally fundamental invariant: the multiplicity. Specifically, we establish a connection between the v-number and the  multiplicity of a homogeneous ideal \(I\) in  \(S = K[x_1, \ldots, x_n]\) over a field \(K\).

\begin{proposition}[Proposition \ref{vNumbervsHilbertMultiplicity}]
    Let $I$ be a zero-dimensional homogeneous ideal in  $S = K[x_1, \dots, x_n]$. Then $\mathrm{v}(I) < e(S/I)$.
\end{proposition}

The paper is organized as follows. In Section \ref{secpreli}, we recall the basic definitions, notation, and preliminary results used throughout the paper. Section \ref{secgradedconjonv} establishes the existence of the asymptotic $\mathrm{v}$-number for Noetherian graded families and identifies it with the asymptotic initial degree. We prove that this limit remains unchanged under integral closure  and provide interpretations through Newton--Okounkov regions for monomial ideals and, under suitable hypotheses, for homogeneous ideals. We also study the corresponding asymptotic invariants of $I^kM$ for finitely generated graded faithful modules $M$ . In Section \ref{section4}, we show that $\operatorname{reg}(I_k)$ and $\mathrm{v}(I_k)$ are eventually  quasi-linear functions of $k$, and derive sufficient conditions for $\mathrm{v}(I_k)<\operatorname{reg}(I_k)$ for large $k$. We then give a formula for a local $\mathrm{v}$-number of stable monomial ideals  and prove that $\mathrm{v}(I)<\operatorname{reg}(I)$ for such ideals . Finally, we establish the inequality $\mathrm{v}(I)<e(S/I)$ for zero-dimensional homogeneous ideals.

\section{Preliminaries}\label{secpreli}

In this section, we provide basic definitions and terminology used in the paper. For unexplained terminology, we refer the interested reader to standard texts \cite{Barvinok02} \cite{HerzogHibi11} \cite{HunekeSwanson06}.  \\ 
Throughout this paper, unless specified stated, we consider $R$  to be a $\mathbb{N}$-graded Noetherian domain with $R_0$ be an infinite field, $S = K[x_1,\ldots,x_n]$ be a standard graded polynomial ring, and $\mathcal{I} = \{I_k\}_{k \geq 0}$ to be a graded family of nonzero ideals.

\begin{definition}
(1) A collection  $\mathcal{I} = \{I_k\}_{k \ge 0}$  of ideals in $R$ is called a \emph{graded family of ideals} if $I_{m}\cdot I_{n}\subseteq I_{m + n}$ for all $m,n\geq 0$.
By convention, we will assume that $I_{0} = R$ when discussing graded families of ideals.
A graded family $\mathcal{I} = \{I_k\}_{k \ge 0}$ is called a \emph{filtration} if $I_{m}\supseteq I_{m + 1}$ for all $m\geq 0$.

\noindent (2) The \emph{Rees algebra} $R(\mathcal{I})$  of a graded family $\mathcal{I} = \{I_k\}_{k \ge 0}$  of ideals in $R$ is defined by
\[
\mathcal{R}(\mathcal{I}) = \bigoplus_{k\geq 0}I_{k}t^{k}\subseteq R[t].
\]

\noindent (3) A graded family $\mathcal{I}$ of ideals in $R$ is said to be \emph{Noetherian} if its Rees algebra $\mathcal{R}(\mathcal{I})$ is a Noetherian ring, i.e., $\mathcal{R}(\mathcal{I})$ is a finitely generated $R$-algebra.

\end{definition}

\begin{theorem}[{\cite[Theorem 2.1]{HerzogHibiTrung07}}]
\label{r}The following conditions are equivalent :
\begin{enumerate}
    \item[(i)] $\mathcal{R}(\mathcal{I})$ is a finitely generated $R$-algebra;
    \item[(ii)] There exists an integer $r \geq 1$ such that $\mathcal{R}(\mathcal{I})^{(r)}= \bigoplus_{k\geq 0}I_{rk}t^{k}$ is a standard graded $R$-algebra;
    \item[(iii)] There exists an integer $r \geq 1$ such that $I_{kr} = I_r^k$ for all $k \geq 0$.
\end{enumerate}

\end{theorem}

\begin{proposition}[{\cite[Proposition 2.8]{HaNguyenNguyen25}}]
\label{Proposition 2.8.}

\noindent Let $\mathcal{R}(\mathcal{I})$ be a finitely generated \(R\)-algebra, and \(r\) is as in Theorem \ref{r}. Then there exist an integer \(N\geq 1\) and ideals \(L_0,\ldots,L_{r-1}\subseteq R\) such that
\[
I_{kr+i} = L_i I_r^{k-N}
\]
for all \(0\leq i\leq r-1\) and all \(k\geq N\) and \(
L_i = I_{Nr + i} \quad \text{for } i = 0,1,\dots, r-1.
\)

In particular, setting \(\ell = Nr\), we have \(I_{k+\ell} = I_k I_\ell\) for all \(k\geq \ell\).

\end{proposition}

\begin{definition}(\cite [Section 2]{HaNguyenNguyen25},\cite[Section 8]{KavehKhovanskii14})
Let $(R,\mathfrak{m})$ be either a local domain or a standard graded domain over a field   ${K}$ and assume that ${K}$ $\simeq R / \mathfrak{m}$.
Let $\operatorname{QF}(R)$ be the quotient field of $R$.
Given a total ordering $>$ in $\mathbb{Z}^{n}$, a valuation $v:\operatorname{QF}(R) \setminus \{0\} \to \mathbb{Z}^{n}$ is a function satisfying the following conditions:
\begin{align*}
(i) &\quad v(fg) = v(f) + v(g) \quad \text{for all } f,g\neq 0; \\
(ii) &\quad v(f + g)\geq \min \{v(f),v(g)\} \quad \text{for all } f,g\neq 0; \\
(iii) &\quad v(\lambda) = 0 \quad \text{for all } 0\neq \lambda \in {K}.
\end{align*}

\end{definition}

A valuation $v:\operatorname{QF}(R) \setminus \{0\}\to \mathbb{Z}^{n}$ is said to have \emph{one-dimensional leaves} if whenever $v(f) = v(g)$, there exists $\lambda \in {K}$ with $v(g + \lambda f) > v(g)$. 
A valuation $v$ with one-dimensional leaves is \emph{good} if the following conditions hold:

\noindent (i) the value semigroup $T = v(R\backslash \{0\})\cup \{0\}$ generates the whole lattice $\mathbb{Z}^{n}$, and its associated cone $C(T)$ is a strongly convex cone;

\noindent (ii) there exists $r_0 > 0$ and a linear function $\ell :\mathbb{R}^{n}\to \mathbb{R}$ such that $C(T)$ lies in $\ell^{-1}(\mathbb{R}_{\geq 0})$ and intersects with $\ell^{-1}(0)$ only at the origin, and for any $f\in R$ if $\ell (v(f))\geq k r_{0}$ for some $k > 0$ then $f\in \mathfrak{m}^{k}$.

\begin{definition}
Let $v$ be a good valuation and let  $\mathcal{I} = \{I_k\}_{k \ge 0}$ be a graded family of ideals in $R$.
For $k\in \mathbb{N}_{>0}$, 
set
\[
\operatorname{val}(I_k) = \{v(f)\mid f\in I_{k}\setminus \{0\} \}.
\]
The \emph{limiting region} and \emph{Newton--Okounkov region} of $\mathcal{I}$ are defined to be
\[
\mathcal{C}(\mathcal{I}) = \bigcup_{k\in \mathbb{N}_{>0}}\mathrm{convex~hull}\left\langle \left\{\frac{\mathbf{a}}{k}\mid \mathbf{a}\in \operatorname{val}(I_k)\right\} \right\rangle \subseteq \mathbb{R}^{n}
\]
and
\[
\Delta (\mathcal{I}) = \overline{\mathcal{C}(\mathcal{I})}\subseteq \mathbb{R}^{n}.
\]

\end{definition}

Now we shall restrict ourselves to a standard graded polynomial ring $S = K[x_1,\ldots,x_n]$. We shall concentrate on good valuations that respect the monomials in $S$. Let $v: \operatorname{QF}(S) \setminus \{0\} \to \mathbb{Z}^n$ be such a valuation. For any monomial $x^a = x_1^{a_1}\cdots x_n^{a_n}$ in $S$ where $a = (a_1,\ldots,a_n) \in \mathbb{Z}_{\geq 0}^n$, we have $v(x^a) = a$. The Gr\"obner valuation of $S$ is an example of a good valuation that respects the monomials. Fix a monomial order on $S$. For a non-zero polynomial $f \in S$, consider all the monomials that appear in $f$ with a non-zero coefficient. Using the monomial order $>$, pick the smallest monomial among them. The value $v(f)$ is simply the exponent vector of the smallest monomial. For example, let $R = K[x,y]$ with the degree reverse lexicographic order where $x > y$. Let $f = 3x^2y + 5xy^2 + 2y^3$. Then $v(f) = (0,3) \in \mathbb{Z}_{\geq 0}^2$. Clearly $v$ can be extended to $\operatorname{QF}(S) \setminus \{0\}$ by $v(f/g) = v(f) - v(g)$.

\begin{definition} 

(\cite [Section 2]{HaNguyenNguyen25}, \cite[Section 4]{KavehKhovanskii14}). Let \(\mathcal{I} = \{I_k\}_{k\geq  0}\) be a graded family of monomial ideals in \(S\). The Newton-Okounkov region of \(\mathcal{I}\) is defined to be

\[
\Delta (\mathcal{I}) = \overline {\bigcup_{k\geq 1}\mathrm{convex~hull}\left\{\frac{\mathbf{a}}{k}\mid x^{\mathbf{a}}\in I_k\right\}} \subseteq \mathbb{R}^n.
\]

\end{definition}

It follows from the definition that \(\Delta (\mathcal{I}) =\overline {\bigcup_{k\in \mathbb{N}}\frac{1}{k}\mathrm{NP}(I_k)}\) . In particular, \(\Delta (\mathcal{I})\) is the closure of the limiting region \(\mathcal{C}(\mathcal{I}) = \bigcup_{k\in \mathbb{N}}\frac{1}{k}\mathrm{NP}(I_k)\) of \(\mathcal{I}\), where        \(\mathrm{NP}(I_k)\) is the Newton polyhedron of $I_k$

\begin{remark}
    For a monomial ideal \(I \subseteq S\):
    \begin{itemize}
        \item If \(\mathcal{I} = \{I^k\}\) is the family of ordinary powers of \(I\),  then \(\Delta(\mathcal{I}) = \operatorname{NP}(I)\).
        \item If \(\mathcal{I} = \{I^{(k)}\}\) is the family of symbolic powers of \(I\), then \(\Delta(\mathcal{I}) = \operatorname{SP}(I)\).
    \end{itemize}
    Thus, both Newton polyhedron $\operatorname{NP}(I)$ and symbolic polyhedron $\operatorname{SP}(I)$ are special cases of the Newton-Okounkov region.
\end{remark}

\begin{definition}

Let \(I\) be an ideal in a ring \(R\). An element \(x \in R\) is called integral over \(I\) if there exists an integer \(m \geq 1\) and elements \(a_i \in I^i\) for \(1 \leq i \leq m\) such that
\[
x^{m} + a_{1}x^{m - 1} + \dots +a_{m - 1}x + a_{m} = 0.
\]
The integral closure of \(I\), denoted by \(\overline{I}\), is an ideal in \(R\), which is the set of all elements in \(R\) that are integral over \(I\) .

\end{definition}

\begin{definition}\label{def:reg-pd-depth}

Let $M$ be a nonzero finitely generated $\mathbb{Z}$-graded module over the standard graded polynomial ring $S=K[x_1,\ldots,x_n]$. Let
\[
\mathbf{F}:\;0\to F_p\to\cdots\to F_1\to F_0\to M\to0
\]
be a minimal graded free resolution of $M$ over $S$, where
\(
F_i=\bigoplus_{j\in\mathbb{Z}}S(-j)^{\beta^S_{i,j}(M)}
\)
for $0\le i\le p$. The \emph{Castelnuovo--Mumford regularity} of $M$, denoted by $\operatorname{reg}(M)$, is given by
\[
\operatorname{reg}(M)=\max\{j-i\mid\beta^S_{i,j}(M)\ne0\};
\]

\noindent Let $R$  be a standard $\mathbb{N}$-graded Noetherian domain. Let $M$ be a finitely generated $\mathbb{Z}$-graded $R$-module. Regularity can also be defined via the vanishing of local cohomology modules with respect to the unique maximal homogeneous ideal
\(
\mathfrak{m} = (x_1,\ldots,x_n).
\)
For every $i\geq 0$, define
\[
a_i(M)
=
\max\{n:H_{\mathfrak{m}}^i(M)_n\neq 0\}.
\]
Then
\[
\operatorname{reg}(M)
=
\max\{a_i(M)+i:i\geq 0\}.
\]

\end{definition}

\section{Asymptotic Behaviour of v-Number of a Noetherian graded family of homogeneous ideals}\label{secgradedconjonv}

In this section, we generalize the results of \cite[Lemma 3.3, Lemma 3.4, Theorem 3.5]{KumarNanduriSaha25} for a Noetherian graded family of homogeneous ideals in $R$.
 We then prove that the asymptotic $\mathrm{v}$-number of $\{I_k\}$ coincides with that of its integral closure $\{\overline{I_k}\}$. Moreover, we provide a combinatorial interpretation of these limits in terms of Newton–Okounkov regions.

\begin{lemma}
Let $\mathcal{I} = \{I_k\}_{k \ge 0}$ be a Noetherian graded family of homogeneous ideals in $R$. Then
\[
\lim_{k \to \infty} \frac{\alpha(I_k)}{k} 
= \inf_{k \ge 1} \frac{\alpha(I_k)}{k}
\]
exists. Moreover, there exists $r \ge 1$ such that
\[
\lim_{k \to \infty} \frac{\alpha(I_k)}{k} = \frac{\alpha(I_r)}{r}.
\]
\end{lemma}

\begin{proof}

The proof follows the same argument as in \cite[Page 7]{VanmathiSarkar25}, where the result was established for a filtration. The same reasoning applies to graded families.

\end{proof}

\begin{lemma}\label{lemma 1.2.}
Let $\mathcal{I} =\{I_k\}_{k \ge 0}$ be a graded family of homogeneous ideals in $R$. 
Then there exists a constant $t \geq 0$ such that for all $k \geq 1$,
\[
\mathrm{v}(I_k) \geq \alpha(I_k) - t.
\]
\end{lemma}

\begin{proof}

 Since $\mathcal{I}$ is a graded family, we have $I_1^k\subseteq I_k$ for every $k\ge1$. If $\mathfrak{p}\in\operatorname{Ass}(I_k)$, then $I_1^k\subseteq I_k\subseteq\mathfrak{p}$. Since $\mathfrak{p}$ is prime, it follows that $I_1\subseteq\mathfrak{p}$, and hence $\alpha(\mathfrak{p})\le\alpha(I_1)$. Set $t=\alpha(I_1)$. Therefore,
\[
c(I_k)=\max\{\alpha(\mathfrak{p})\mid\mathfrak{p}\in\operatorname{Ass}(I_k)\}\le t.
\]
By \cite[Lemma 3.1]{BiswasMandalSaha24}, we obtain
\[
\mathrm{v}(I_k)\ge\alpha(I_k)-c(I_k)\ge\alpha(I_k)-t,
\]
as required.
\end{proof}

\begin{theorem}\label{Proposition 1.4}
Let \(\mathcal{I} = \{I_k\}_{k \geq 0}\) be a Noetherian graded family of homogeneous ideals in  \(R\). Then
\[
\lim_{k \to \infty} \frac{\operatorname{v}(I_k)}{k}
\]
exists, and there exists \(r \geq 1\) such that
\[
\lim_{k \to \infty} \frac{\operatorname{v}(I_k)}{k} = \frac{\alpha(I_r)}{r}.\]
\end{theorem}

\begin{proof}
Since \(\mathcal{I}\) is Noetherian, there exists an integer \(r \ge 1\) such that \(I_{nr} = I_r^n\) for all \(n \ge 0\). For any integer \(k \geq 0\), we can write uniquely \(k = rn + i\) with \(0 \leq i \leq r-1\). Using Proposition \ref{Proposition 2.8.}, there exists a sufficiently large integer \(N\) such that for all \(0 \leq i \leq r-1\) and all \(n > N\),
\[
I_{rn+i} = I_{Nr+i}I_r^{n-N}.
\]
Write \(L_i = I_{Nr+i}\) for all \(0 \leq i \leq r-1\). Fix a homogeneous element \(f \in I_r\) with \(\deg f = \alpha(I_r)\), then \(f^n \in I_{rn}\) and \(\alpha(I_{rn}) = \alpha(I_r^n) = n\alpha(I_r)\). Moreover, \(L_i f^{n-N} \subseteq I_{rn+i}\) and \(f^{n-N-1} \notin I_{rn+i}\) for all \(n > N\). Applying \cite[Proposition 4.3]{{FicarraSgroi23}} we get,
\[\operatorname{v}(I_{rn+i}) \leq \operatorname{v}\big(I_{rn+i} : f^{n-N-1}\big) + (n-N-1)\alpha(I_r). \tag{1}\]
Now, for a fixed \(i \geq 0\), we have the following increasing sequence of ideals in \(R\)
\[
(I_{rn+i} : f^{n-N-1}) \subseteq (I_{r(n+1)+i} : f^{n-N}) \subseteq (I_{r(n+2)+i} : f^{n-N+1}) \subseteq \cdots.
\]
Since \(R\) is Noetherian, there exists a positive integer \(p_i\) such that
\[
(I_{rn+i} : f^{n-N-1}) = (I_{rp_i+i} : f^{p_i-N-1}) \quad \text{for all } n \geq p_i,
\]
and a constant \(c_i\) with \(\operatorname{v}(I_{rn+i} : f^{n-N-1}) = c_i\) for \(n \geq p_i\). Hence from equation ({1}) we get
\[
\operatorname{v}(I_{rn+i}) \leq c_i + (n-N-1)\alpha(I_r) \quad (n \geq p_i).\tag{2}
\]
Using Lemma \ref{lemma 1.2.},  Equation (2) and the equality \(\alpha(I_{rn+i}) = \alpha(L_i) + (n-N)\alpha(I_r)\) we obtain for \(k = rn + i\) with \(n \geq p_i\):
\[
\frac{\alpha(L_i) + (n-N)\alpha(I_r) - t}{rn+i} \leq \frac{\operatorname{v}(I_{rn+i})}{rn+i} \leq \frac{c_i + (n-N-1)\alpha(I_r)}{rn+i}.
\]
As \(n \to \infty\), we have \(k = rn + i \to \infty\) and both bounds converge to \(\frac{\alpha(I_r)}{r}\). Since this holds for each residue class \(i = 0,1,\ldots,r-1\), we conclude that
\[
\lim_{k \to \infty} \frac{\operatorname{v}(I_k)}{k} = \frac{\alpha(I_r)}{r},
\]
which completes the proof.
\end{proof}

The following example demonstrates Theorem \ref{Proposition 1.4} where we do not consider our family to be a filtration.

\begin{example}
 Let  
$I = (x_1,\dots,x_n) = \mathfrak{m}$, where $n\geq2$ be the unique homogeneous maximal ideal in $S$ and $J = (x_1)$, 
$I_0 = S$, 
\[
I_k = 
\begin{cases}
I^k, & \text{if } k \text{ is even} \\
JI^{k}, & \text{if } k \text{ is odd}
\end{cases}
\]
is a Noetherian graded family of homogeneous ideals in $S$.\\
When $k$ is even: 
$I_k = \mathfrak{m}^k$, $\mathrm{Ass}(I_k) = \mathfrak{m}$, 
$I_k : \mathfrak{m} = \mathfrak{m}^{k-1}$, 
$x_1^{k-1} \in (I_k : \mathfrak{m}) \setminus I_k$ has degree $k - 1$. 
Moreover by \cite[Lemma 1.2]{Conca24} $\mathrm{v}(I_k) = \mathrm{v}_{\mathfrak{m}}(I_k) = \min \{ w : \left( (I_k : \mathfrak{m}) / I_k \right)_w \neq 0 \} = k-1$.\\
When $k$ is odd: 
$I_k = (x_1)\mathfrak{m}^{k} = (x_1) \cap \mathfrak{m}^{k+1}$, 
$\mathrm{Ass}(I_k) = \{ \mathfrak{p} = (x_1), \mathfrak{m} \}$, 
$I_k : \mathfrak{m} = ((x_1)\mathfrak{m}^{k}) : \mathfrak{m} = (x_1)\mathfrak{m}^{k-1}$, 
$x_1^{k} \in (I_k : \mathfrak{m}) \setminus I_k$ is of degree $k$, 
$\mathrm{v}_{\mathfrak{m}}(I_k) = k$, 
$I_k : \mathfrak{p} = ((x_1)\mathfrak{m}^{k}) : (x_1) = \mathfrak{m}^{k}$, 
$I_k : \mathfrak{m}^{\infty} = (x_1)$, 
$I_k : (\mathfrak{p} + \mathfrak{m}^{\infty}) = (I_k : \mathfrak{p}) \cap (I_k : \mathfrak{m}^{\infty}) = \mathfrak{m}^{k} \cap (x_1)$, 
$x_2^{k} \in \mathfrak{m}^{k} \setminus (\mathfrak{m}^{k} \cap (x_1))$, 
$\mathrm{v}_{\mathfrak{p}}(I_k) = \min \{ w: \left((I_k:\mathfrak{p}\right) / I_k: (\mathfrak{p} + \mathfrak{m}^{\infty}))_w \neq 0 \} = k$, 
$\mathrm{v}(I_k) = \min \{ \mathrm{v}_{\mathfrak{m}}(I_k), \mathrm{v}_{\mathfrak{p}}(I_k) \} = k$.
Therefore 
\[
\mathrm{v}(I_k) = 
\begin{cases}
k-1, & \text{when } k \text{ is even} \\
k, & \text{when } k \text{ is odd}
\end{cases}
\quad\text{and}\quad
\alpha(I_k) = 
\begin{cases}
k, & \text{when } k \text{ is even} \\
k+1, & \text{when } k \text{ is odd}
\end{cases}
\]
and 
\[
\lim_{k \to \infty} \frac{\mathrm{v}(I_k)}{k} = \lim_{k \to \infty} \frac{\alpha(I_k)}{k} = 1.
\]
\end{example}

In \cite[Theorem 4.1]{HaNguyenNguyen25}, H. T. Hà et al. showed that for a graded family $\mathcal{I} = \{I_k\}_{k \ge 0}$ of $\mathfrak{m}$-primary homogeneous ideals,
\(
\lim\limits_{k \to \infty} \frac{\operatorname{reg}(I_k)}{k} = \lim\limits_{k \to \infty} \frac{d(I_k)}{k}.
\)
However, an analogous result does not hold for the $\mathrm{v}$-number. We illustrate this with an example.

\begin{example}
    
    Let \(I_0 = S\), \(I_1 = (x,y^2)\), and \(I_k = (x^k, y)\) for \(k \ge 2\) be a non-Noetherian graded family of \(\mathfrak{m}\)-primary homogeneous ideals in \(S = K[x,y]\). 
    Note that \(\operatorname{Ass}(I_k) = (x,y) = \mathfrak{m}\), so \(\mathrm{v}(I_k) = \mathrm{v}_{\mathfrak{m}}(I_k)\). 
    Moreover, \(I_k : (x,y) = (x^{k-1}, y)\), and \(x^{k-1} \in (I_k : \mathfrak{m}) \setminus I_k\) has degree \(k-1\).
    By \cite[Lemma 1.2]{Conca24}, 
    \[
    \mathrm{v}(I_k) = \min \{ w : \left((I_k : \mathfrak{m}\right) / I_k)_w \neq 0 \} = k-1,
    \quad
    \lim_{k \to \infty} \frac{\mathrm{v}(I_k)}{k} = 1 \neq \lim_{k \to \infty} \frac{\alpha(I_k)}{k} = 0.
    \]
\end{example}

\noindent Next, we prove the following lemma, which is essential for establishing our main result.

\begin{lemma}\label{alpha(I)=alpha(I^bar)}

Let $I$ be a non-zero homogeneous ideal in $R$. Then $\alpha(I) = \alpha(\overline{I}).$
\end{lemma}

\begin{proof}
    
Since \(I\) is homogeneous, its integral closure \(\overline{I}\) is also homogeneous by \cite[Corollary 5.2.2]{HunekeSwanson06}. Again $I \subseteq \overline{I}$, so $\alpha(\overline{I}) \leq \alpha(I)$.

Conversely, let $f \in \overline{I}$ be a non-zero homogeneous element of degree $d$ such that $\deg(f) = \alpha(\overline{I}) = d$. Then 
\[
f^n + c_1 f^{n-1} + c_2 f^{n-2} + \cdots + c_n = 0
\]
where $c_i \in I^i$ and $c_i = \sum_j c_{i,j}$ where $c_{i,j}$ $\in I^i$ is a homogeneous element of degree $j$. Thus
\[
f^n + \sum_{i=1}^n \left(\sum_j c_{i,j}\right) f^{n-i} = 0.
\]
Therefore, 
\[
f^n + a_1 f^{n-1} + a_2 f^{n-2} + \cdots + a_n = 0
\]
is a homogeneous equation of degree $nd$ where $a_i = c_{i,id}$  $\in I^i$ for all $1 \leq i \leq n$, and for some $a_j \neq 0$, we have deg $a_j = jd$. It implies $\alpha(\overline{I}) \geq \alpha(I)$.
\end{proof}

\begin{theorem}\label{mainresultofhomogeneousideals}
Let $\mathcal{I}=\{I_k\}_{k \ge 0}$ be a Noetherian graded family of homogeneous ideals in  $R$. Then
\[
\lim_{k\to\infty}\frac{\mathrm{v}(\overline{I_k})}{k}
   =\lim_{k\to\infty}\frac{\alpha(\overline{I_k})}{k}
   =\lim_{k\to\infty}\frac{\mathrm{v}(I_k)}{k}
   =\lim_{k\to\infty}\frac{\alpha(I_k)}{k}.
\]
\end{theorem}

\begin{proof}
The family $\overline{\mathcal{I}}=\{\overline{I_k}\}_{k \ge 0}$ is again Noetherian and  graded by \cite[Lemma 2.10]{{HaNguyenNguyen25}}   and \cite[Corollary 6.8.6]{{HunekeSwanson06}}.  
By  Theorem \ref{r}, there exists $r\ge 1$ such that $\overline{I_{rk}}=(\overline{I_r})^{\,k}$ for all $k\ge 0$.  
Applying Theorem \ref{Proposition 1.4} to $\overline{\mathcal{I}}$ yields
\[
\lim_{k\to\infty}\frac{\mathrm{v}(\overline{I_k})}{k}=\displaystyle\lim_{k\to\infty}\frac{\alpha(\overline{I_k})}{k}=\frac{\alpha(\overline{I_r})}{r}.
\]
By Lemma \ref{alpha(I)=alpha(I^bar)}, $\alpha(\overline{I_k})=\alpha(I_k)$ for every $k$. Combining these equalities completes the proof.
\end{proof}

\begin{theorem}\label{analogue}
Let $\mathcal{I}=\{I_k\}_{k \ge 0}$ be a Noetherian graded family of monomial ideals in $S $. Then
\[
\lim_{k\to\infty}\frac{\mathrm{v}(\overline{I_k})}{k}
   =\lim_{k\to\infty}\frac{\alpha(\overline{I_k})}{k}
   =\lim_{k\to\infty}\frac{\mathrm{v}(I_k)}{k}
   =\lim_{k\to\infty}\frac{\alpha(I_k)}{k}
   =\lambda(\Delta(\mathcal{I})),
\]
where $\lambda(\Delta(\mathcal{I})):=\min\{|\mathbf{v}|:\mathbf{v}\text{ is a vertex of }\Delta(\mathcal{I})\}$.
\end{theorem}

\begin{proof}
Since \(\mathcal{I}\) is Noetherian, there exists a positive integer $r$ such that the $r$-th Veronese subalgebra of $\mathcal{R}(\mathcal{I})$ is standard graded. By \cite[Theorem 4.5]{CDF+}, we have
\[
\lim_{k \to \infty} \frac{\alpha(I_k)}{k} = \inf_{k \ge 1} \frac{\alpha(I_k)}{k} = \inf_{k \ge 1} \left\{ \min\{y_1 + \cdots + y_n \mid (y_1,\dots,y_n) \in \tfrac{1}{k}\operatorname{NP}(I_k)\} \right\}
\]
and
\[
\lim_{k \to \infty} \frac{\alpha(I_k)}{k} = \lambda(\Delta(\mathcal{I})).
\]

\end{proof}
Moreover, by \cite[Theorem 3.1]{HaNguyen}, we can replace $\Delta(\mathcal{I})$ with $\frac{1}{r}\operatorname{NP}(I_r)$. In the particular case where $\mathcal{I} = \{I^n\}$, the limit $\lambda(\Delta(\mathcal{I}))$ is an non-negative integer. In general, \(\lambda(\Delta(\mathcal{I}))\) need not be equal to a non‑negative integer.

\begin{example}
    
     Let \(I_k = (x^{a_k}) \subseteq K[x]\) where 
    \[
    a_k = \begin{cases} 
    \frac{k}{2}, & \text{if } k \text{ is even}, \\
    \frac{k+3}{2}, & \text{if } k \text{ is odd},
    \end{cases}
    \]
    be a Noetherian graded family of homogeneous ideals in \(K[x]\). Then
    \[
    \mathrm{v}(I_k) = \begin{cases} 
    \frac{k}{2} - 1, & \text{if } k \text{ is even}, \\
    \frac{k+1}{2}, & \text{if } k \text{ is odd},
    \end{cases}
    \qquad
    \alpha(I_k) = \begin{cases} 
    \frac{k}{2}, & \text{if } k \text{ is even}, \\
    \frac{k+3}{2}, & \text{if } k \text{ is odd}.
    \end{cases}
    \]
    Consequently, \(\Delta(\mathcal{I}) = [1/2, \infty)\) with vertex \(1/2\).  
    Therefore, \(\lambda(\Delta(\mathcal{I})) = 1/2\).
\end{example}

In \cite{HaNguyenNguyen25}, H. T. Hà et al. related asymptotic regularity to the Newton–Okounkov region for a Noetherian graded family of monomial ideals and \(\mathfrak{m}\)-primary homogeneous ideals. Motivated by their result, we prove the following:

\begin{theorem}\label{NORforvaluations}

Let $\mathcal{I} = \{I_k\}_{k \ge 0}$ be a Noetherian graded family of homogeneous ideals, 
and let $v$ be a good valuation that respects the monomials in $S$. 
Let $\Delta(\mathcal{I})$ be the Newton--Okounkov region of 
$\mathcal{I}$ defined by $v$. Suppose that the limiting region $\mathcal{C}(\mathcal{I})$ is a polyhedron. Then
\[
\lim_{k \to \infty} \frac{\mathrm{v}(\overline{I_k})}{k} 
= \lim_{k \to \infty} \frac{\alpha(\overline{I_k})}{k}
= \lim_{k \to \infty} \frac{\mathrm{v}(I_k)}{k}
= \lim_{k \to \infty} \frac{\alpha(I_k)}{k}
= \lambda(\Delta(\mathcal{I})),
\]
where $\lambda(\Delta(\mathcal{I})) := \min\{\,|\mathbf{v}| : \mathbf{v} \text{ is a vertex of } \Delta(\mathcal{I})\,\}$.

\end{theorem}

\begin{proof}
 Let $v$ be a good valuation that respects the monomials in $S$, and define
\[
\widetilde{I_k} = \big\langle \{x^{v(f)} \mid f \text{ is a homogeneous element of } I_k  \backslash   \{0\}\} \big\rangle,
\]
Note that  $\widetilde{I_k}$ is a monomial ideal. As valuation $v$ respects the monomials in S, we have
$\operatorname{val}(I_k) = \operatorname{val}(\widetilde{I_k}) = \operatorname{Ex}(\widetilde{I_k})$, where \(\operatorname{Ex}(\widetilde{I_k})\) denotes the set of exponent vectors of monomials in \(\widetilde{I_k}\).
This implies that $\alpha(I_k) = \alpha(\widetilde{I_k})$, $\mathcal{C}(\mathcal{I}) = \mathcal{C}(\widetilde{\mathcal{I}})$, and $\Delta(\mathcal{I}) = \Delta(\widetilde{\mathcal{I}})$ where $\widetilde{\mathcal{I}} = \{\widetilde{I_k}\}_{k \geq 0}$ is a graded family of monomial ideals in $S$. Consequently,
\[
\lim_{k \to \infty} \frac{\alpha(I_k)}{k} = \lim_{k \to \infty} \frac{\alpha(\widetilde{I_k})}{k}.
\]

\noindent Since $\mathcal{C}(\widetilde{\mathcal{I}})$ is a polyhedron, by \cite[Theorem 3.4]{{HaNguyen}} \(\overline{\mathcal{I}}\) is Noetherian and by Proposition \ref{analogue}, we have
\[ \lim_{k \to \infty} \frac{\mathrm{v}(I_k)}{k}
= \lim_{k \to \infty} \frac{\alpha(I_k)}{k}
= 
\lim_{k \to \infty} \frac{\alpha(\widetilde{I_k})}{k} = \lambda\bigl(\Delta(\widetilde{\mathcal{I}})\bigr)=\lambda\bigl(\Delta({\mathcal{I}})\bigr)
\]
which completes the proof.

\end{proof}

Before proceeding to the next proposition, we recall some necessary preliminaries from \cite[Section 1]{{FiorindoGhosh}}. Let $R$ be a Noetherian $\mathbb{N}$-graded domain. Let $M$ be a finitely generated $\mathbb{Z}$-graded $R$-module and $I$ be a homogeneous ideal of $R$. Let $N$ be a graded submodule of $M$. Let $J$ be a reduction ideal of $I$, i.e., there exists a non-negative integer $n$ such that $JI^n = I^{n+1}$. It implies $\alpha(I) = \alpha(J)$. Let $J$ be generated by homogeneous elements $y_1, \dots, y_c$ of degree $d_1 \le \dots \le d_c$ respectively. The set of associated prime ideals of the $R$-module $M$ is denoted by $\operatorname{Ass}_R(M)$. For each $\mathfrak{p} \in \operatorname{Ass}_R(M)$, the $\mathrm{v}$-number of $M$ at $\mathfrak{p}$ is defined as
\[
\mathrm{v}_{\mathfrak{p}}(M) := \inf \{ u \mid \text{there exists } x \in M_u \text{ such that } \mathfrak{p} = (0 :_R x) \}.
\]
Then, the $\mathrm{v}$-number of $M$ is
\[
\mathrm{v}(M) := \inf \{ \mathrm{v}_{\mathfrak{p}}(M) \mid \mathfrak{p} \in \operatorname{Ass}_R(M) \} \quad \text{and} \quad \alpha(M) := \inf \{ t \mid M_t \neq 0 \}.
\]

\begin{proposition}\label{vnumberofmodules}
Let $I$ be a homogeneous ideal in $R$. 
Then, for any finitely generated graded faithful $R$-module $M$, we have
\[
\lim_{k \to \infty} \frac{\mathrm{v}(I^k M)}{k} 
= \lim_{k \to \infty} \frac{\alpha(I^k M)}{k} 
= \lim_{k \to \infty} \frac{\mathrm{v}(I^k)}{k} 
= \alpha(I)
\]
\end{proposition}

\begin{proof}

Let \(J\) be a homogeneous reduction ideal  of \(I\). Write \(J = (y_1,\dots,y_c)\) with each \(y_j\) homogeneous and \(\deg(y_j) = d_j\), ordered so that \(d_1 \le d_2 \le \cdots \le d_c\). Since \(\alpha(I) = \alpha(J)\), we have \(\alpha(I) = d_1\).

Consider the graded module \(H = \bigoplus_{k\ge 0} I^kM\). Since \(R\) is a domain and \(M\) is faithful, \(I^kM \neq 0\) for all \(k\); hence the set \(\operatorname{Ass}_R(I^kM)\) is nonempty. By \cite[Theorem 2.11]{{FiorindoGhosh}}, both functions \(\operatorname{\alpha}(I^kM)\) and \(\mathrm{v}(I^kM)\) are eventually linear in \(k\), and they share the same leading coefficient \(d_\delta\), where  
\[
\delta = \inf\{\,j : y_j \notin \sqrt{\operatorname{ann}_{\mathcal{R}(J)}(H)}\,\}.
\]
We claim that \(\delta = 1\). Suppose that some power of \(y_1\) annihilates \(H\); i.e., \(y_1^t H = 0\) for some \(t > 0\). Then in particular \(y_1^t M = 0\). Since \(M\) is faithful, this forces \(y_1^t = 0\) in \(R\), contradicting the fact that \(R\) is a domain. Hence \(y_1 \notin \sqrt{\operatorname{ann}_{\mathcal{R}(J)}(H)}\), so \(\delta = 1\) and consequently \(d_\delta = d_1 = \alpha(I)\). Thus we get\\
\begin{equation}\label{eq1}
    \lim_{k \to \infty} \frac{\mathrm{v}(I^k M)}{k} 
= \lim_{k \to \infty} \frac{\alpha(I^k M)}{k} 
= \alpha(I)
\end{equation}
Now combining (\ref{eq1}) and applying Theorem \ref{mainresultofhomogeneousideals}, we have 
\[
\lim_{k \to \infty} \frac{\mathrm{v}(I^k M)}{k} 
= \lim_{k \to \infty} \frac{\alpha(I^k M)}{k} \lim_{k \to \infty} \frac{\mathrm{v}(I^k)}{k} = \alpha(I)
.\]

\end{proof}

\section{A comparison of v-number with regularity and Multiplicity}\label{section4}

\noindent In this section we prove that  the $\mathrm{v}$-number of a Noetherian graded family of homogeneous ideals in $R$ is a quasi-linear function. Then we compare the \(\mathrm{v}\)-number with regularity and multiplicity for some classes of ideals.

\begin{definition}

A numerical function $f:\mathbb{Z}_{\geq 0}\to \mathbb{Z}_{\geq 0}$ is called a \emph{quasi-linear function} (of period $r$) if there exist a positive integer $r$ and linear functions $f_{i}(k) = a_{i}k + b_{i}$, $i = 0,\ldots ,r - 1$ such that $f(k) = f_{i}(k)$ if $k\equiv i$ (mod $r$).

\end{definition}

Let $I$ be a homogeneous ideal in an $\mathbb{N}$-graded domain $R$. Let $M$ be a finitely generated graded $R$-module. A homogeneous ideal $J \subseteq I$ is an $M$-\emph{reduction} of $I$ if $I^{n+1}M = JI^nM$ for some $n \geq 0$. Define
\[\rho_M(I) = \min\{ d(J) \mid J \text{ is an } M\text{-reduction of } I\}\]
where for a graded module $N$, $d(N)$ denotes the maximal degree of a minimal homogeneous generating set of $N$.

\begin{theorem}\label{reg/vnumberisquasilinear}
   Assume, in addition, that $R$ is standard graded. Let $\mathcal{I} = \{ I_k \}_{k \geq 0}$ be a Noetherian graded family of homogeneous ideals in $R$. Then
     $\operatorname{reg}(I_k)$ is eventually a quasi-linear function in $k$.
\end{theorem}

\begin{proof}

Since $\mathcal{I}$ is a Noetherian graded family, its Rees algebra $\mathcal{R}(\mathcal{I})$ is finitely generated over $R$. By Proposition \ref{Proposition 2.8.}, there exist positive integers $r$ and $N$, and ideals $L_0, L_1, \dots, L_{r-1} \subseteq R$ such that for all $m \geq N$ and $0 \leq i \leq r-1$, we have
\[
I_{mr+i} = I_r^{m-N} L_i,
\]
where $L_i = I_{Nr+i}$.

\noindent For each fixed $i$ with $0 \leq i \leq r-1$, consider the module $M = L_i$, which is a finitely generated graded faithful $R$-module. By the proof of \cite[Theorem 3.3]{HaNguyenNguyen25}, we have
\[
\lim_{n \to \infty} \frac{\operatorname{reg}(I_r^n L_i)}{n} = \lim_{n \to \infty} \frac{\operatorname{reg}(I_r^n)}{n} = \rho_M(I_r) = \rho_R(I_r) = \rho.
\]
for every $i$. Moreover, by  \cite[Theorem 3.2]{{TrungWang}}, for each $i$ there exists an integer $e_i$ such that for all sufficiently large $n$,
\[
\operatorname{reg}(I_r^nL_i) = \rho n + e_i.
\]

\noindent Now write $k = mr + i$ with $0 \leq i \leq r-1$ and $m \gg 0$. Then $I_k =  I_r^{m-N}L_i$ for $m \geq N$. Applying the above with $n = m-N$, we obtain for all $m$ sufficiently large,
\[
\operatorname{reg}(I_k) = \operatorname{reg}(I_r^{m-N}L_i) = \rho(m-N) + e_i = \rho m + (e_i - \rho N).
\]

\noindent Since $m = \frac{k-i}{r}$, substituting gives
\[
\operatorname{reg}(I_k) = \rho\left(\frac{k-i}{r}\right) + (e_i - \rho N) = \frac{\rho}{r} k + \left(e_i - \rho N - \frac{\rho i}{r}\right).
\]

\noindent Let $A = \frac{\rho}{r}$ and $B_i = e_i - \rho N - \frac{\rho i}{r}$. Then for all $k \gg 0$ with $k \equiv i \pmod{r}$, we have
\[
\operatorname{reg}(I_k) = A k + B_i.
\]

\noindent Thus $\operatorname{reg}(I_k)$ is eventually a quasi-linear function in $k$, with the coefficient $A$ independent of the residue class modulo $r$. This completes the proof.

\end{proof}

 Note \cite[Lemma 1.1]{FicarraSgroi24}, states that if $\mathcal{I} = \{ I_k \}_{k \geq 0}$ be a Noetherian filtration of homogeneous ideals in $R$. Then there exists an integer $r > 0$ such that $\operatorname{Ass}(I_{rk+i}) = \operatorname{Ass}(I_{r(k+1)+i})$ for all $i = 0,\ldots,r-1$ and $k \gg 0$. The same proof follows for a Noetherian graded family of homogeneous ideals in R.  Again, \cite[Corollary 4.2]{Ratliff76} the Ratliff property states that if $I$ is an ideal of a Noetherian graded domain $R$, then we have $(I^{k+1}: I) = I^k$ for all $k \gg 0$.

\begin{lemma}\label{pI^k+1:p}
Let \(I \subseteq S\) be a monomial ideal, and let \(\mathfrak{p} \subseteq S\) be a nonzero monomial prime ideal. Then, for every integer \(k \ge 0\),
\[
(\mathfrak{p} I^{k+1} : \mathfrak{p}) \subseteq I^k.
\]
\end{lemma}

\begin{proof}
Since all ideals involved are monomial, it suffices to check monomial elements. Let \(u\) be a monomial such that \(u \in (\mathfrak{p} I^{k+1} : \mathfrak{p})\). Pick any variable \(x_i \in \mathfrak{p}\). Then \(u x_i \in \mathfrak{p} I^{k+1}\), so there exist \(x_j \in \mathfrak{p}\) and a monomial \(m \in I^{k+1}\) such that
\[
u x_i = x_j m. \tag{1}
\]
If \(j = i\), then \(u = m \in I^{k+1} \subseteq I^k\), and we are done.

 Assume now that \(j \ne i\). Then \(x_i\) divides \(m\). Write \(m = m_1 m_2 \cdots m_{k+1}\), where each \(m_\ell\) is a monomial belonging to \(I\). Since \(x_i \mid m\), after relabelling we may write \(m_1 = x_i m_1'\) for some monomial \(m_1'\). From (1) we obtain
\[
u = \frac{x_j}{x_i} m = x_j m_1' m_2 \cdots m_{k+1}.
\]
The product \(m_2 \cdots m_{k+1}\) is a product of \(k\) generators of \(I\), hence \(u \in I^k\). This proves the lemma. \qedhere
\end{proof}

\begin{theorem}\label{visquasilinear}

Let \(\mathcal{I} = \{I_k\}_{k \geq 0}\) be a Noetherian graded family of monomial ideals in \(S \). Then \(\mathrm{v}(I_k)\) is  eventually  a quasi-linear function in \(k\). 
\end{theorem}

\begin{proof}
Let \(r > 0\) be a positive integer such that \(\mathcal{R} = \mathcal{R}(\mathcal{I})^{(r)} = \bigoplus_{k \geq 0} I_{rk}\) is a standard graded \(S\)-algebra. For each \(i \in \{0, 1, \dots, r-1\}\), let \(\operatorname{Ass}^{\infty}(\mathcal{I})_i\) be the set of those prime ideals \(\mathfrak{p} \in \operatorname{Spec}(S)\) such that \(\mathfrak{p} \in \operatorname{Ass}(I_{rk + i})\) for all \(k \gg 0\). Hence, for each \(i \in \{0, \dots, r-1\}\) and all \(k \gg 0\), we have
\[
\mathrm{v}(I_{rk + i}) = \min_{\mathfrak{p} \in \operatorname{Ass}^{\infty}(\mathcal{I})_i} \mathrm{v}_{\mathfrak{p}}(I_{rk + i}).
\]

\noindent Now fix \(i \in \{0, \dots, r-1\}\). Let \(\mathfrak{p} \in \operatorname{Ass}^{\infty}(\mathcal{I})_i\), 
\(
X_{\mathfrak{p}} = \{ \mathfrak{p}_1 \in \operatorname{Ass}^{\infty}(\mathcal{I})_i : \mathfrak{p} \subsetneq \mathfrak{p}_1 \}.
\)
Let \(\mathfrak{q}\) be the product of the elements in \(X_{\mathfrak{p}}\), with the convention that \(\mathfrak{q} = S\) if \(X_{\mathfrak{p}}\) is empty.

\noindent Let \(\mathcal{R}' = I_{r(N+1)+i} \mathcal{R} \) , by Proposition \ref{Proposition 2.8.} we have \(\mathcal{R}' = \bigoplus_{k \geq 0} I_{r(k+N+1)+i}\), \(\mathfrak{p} \mathcal{R} = \bigoplus_{k \geq 0} \mathfrak{p} I_{rk}\), and \(\mathfrak{q} \mathcal{R} = \bigoplus_{k \geq 0} \mathfrak{q} I_{rk}\). Now set
\[
A = (\mathcal{R}' :_{\mathcal{R}} \mathfrak{p} \mathcal{R}) \quad \text{and} \quad B = \mathcal{R}' :_{\mathcal{R}} (\mathfrak{p} + \mathfrak{q}^{\infty}) \mathcal{R}
\]
such that \(H = A / B\) is a finitely generated bigraded \(\mathcal{R}\)-module. Then
\[
H_{(*,k)} = \frac{(I_{r(k+N+1)+i} : \mathfrak{p}) \cap I_{rk}}{(I_{r(k+N+1)+i} : \mathfrak{p} + \mathfrak{q}^{\infty}) \cap I_{rk}}.
\]

Since \(\mathfrak{p} \in \operatorname{Ass}(I_{rk+i})\) for all \(k \gg 0\), we have \(I_{rk+i} \subseteq \mathfrak{p}\). Then \(I_r^k I_i \subseteq \mathfrak{p}\).

\noindent \textbf{Case 1:} \(I_r \subseteq \mathfrak{p}\). Let \(x \in (I_{r(k+N+1)+i} : \mathfrak{p})\). This implies \(x I_r \subseteq I_{r(k+N+1)+i} = I_{rN+i} I_r^{k+1}\) and consequently \(x I_r \subseteq I_r^{k+1}\). Hence \(x \in (I_r^{k+1} : I_r) = I_r^k\).

\noindent \textbf{Case 2:} \(I_i \subseteq \mathfrak{p}\) for \(i = 1, \dots, r-1\). Since \(\mathfrak{p} \in \operatorname{Ass}(I_{rk+i})\) for all \(k \gg 0\), we have \(I_{rN+i} I_r^{k+1}\subseteq \mathfrak{p}\) such that \(I_{Nr+i} \subseteq \mathfrak{p}\). This implies
\(
(I_{r(k+N+1)+i} : \mathfrak{p}) = (I_{rN+i} I_r^{k+1} : \mathfrak{p}) \subseteq (\mathfrak{p} I_r^{k+1} : \mathfrak{p}),
\)
and by Lemma \ref{pI^k+1:p}, \((\mathfrak{p}I_r^{k+1} : \mathfrak{p}) \subseteq I_r^{k+1} \subseteq I_r^k\).

\noindent Therefore,
\[
H_{(*,k)} = \frac{(I_{r(k+N+1)+i} : \mathfrak{p})}{(I_{r(k+N+1)+i} : \mathfrak{p} + \mathfrak{q}^{\infty}) }.
\]
 for all \(k \gg 0\). By applying \cite[Lemma 1.2]{Conca24} , we get
\[
\mathrm{v}_{\mathfrak{p}}(I_{rk+i}) = \alpha(H_{(*,k-N-1)}).
\]
A straightforward modification of the arguments used in the proofs of
\cite[Lemma 3.3]{BrunsConcaVarbaro} and \cite[Theorem 8.3.4]{BrunsConcaRaicuVarbaro2022} shows that, $\mathrm{v}_{\mathfrak{p}}(I_{rk + i})$
is eventually linear in \(k\).

\end{proof}

\noindent We now extend Theorem~\ref{visquasilinear} to Noetherian graded families of homogeneous ideals; for this purpose, we first prove the following lemma.

\begin{lemma}\label{ArtinReeslemma}
Let $R$ be a Noetherian ring, $I\subseteq R$ an ideal, and let $a\in R$ be a non-zero divisor. Then there exists an integer $s\ge 1$ such that, for every $k\ge 0$,
\[
(I^{k+s}:_R a)\subseteq I^k.
\]
\end{lemma}

\begin{proof}
By the Artin-Rees Lemma, there exists an integer $s\ge 1$ such that
\(
I^n\cap aR = I^{n-s}(I^s\cap aR)
\)
for all $n\ge s$. Set $n=k+s$. Then
\(
I^{k+s}\cap aR = I^k(I^s\cap aR)
\)
for all $k\ge 0$.

\noindent For every integer $m\ge 0$, we have
\(
I^m\cap aR = a(I^m:_R a).
\)
Now,
\[
a(I^{k+s}:_R a) = I^{k+s}\cap aR = I^k(I^s\cap aR) = I^k\cdot a(I^s:_R a) \subseteq aI^k.
\]
Since $a$ is a non-zero divisor, it follows that
\(
(I^{k+s}:_R a)\subseteq I^k.
\)
\end{proof}

\noindent In the next theorem, we adopt the notation and setup of the proof of Theorem~\ref{visquasilinear}, with \(R\) in place of \(S\).

\begin{theorem}\label{vnumberisquasilinearforhomogeneous}
Let \(\mathcal{I} = \{I_k\}_{k \geq 0}\) be a Noetherian graded family of homogeneous ideals in  \(R\). Then \(\mathrm{v}(I_k)\) is eventually a quasi-linear function in \(k\). 
\end{theorem}

\begin{proof}

Since $\mathfrak{p}\in\operatorname{Ass}(I_{rk+i})$ for all $k\gg 0$, we have $I_{rk+i}\subseteq\mathfrak{p}$. Since $I_{rk+i}\neq 0$, we have $\mathfrak{p}\neq (0)$. Let $0\neq a\in\mathfrak{p}$. Since $R$ is a domain, $a$ is a non-zero divisor. By Lemma \ref{ArtinReeslemma}, there exists an integer $s\geq 1$ such that, for every $k\geq 0$,
\[
(I_r^{k+s}:_R a)\subseteq I_r^k.
\]

\noindent Let \(\mathcal{R}' = I_{r(N+s+1)+i} \mathcal{R}\). By Proposition \ref{Proposition 2.8.}, we have  
\( 
\mathcal{R}' = \bigoplus_{k \geq 0} I_{r(k+N+s+1)+i}, \quad  
\mathfrak{p} \mathcal{R} = \bigoplus_{k \geq 0} \mathfrak{p} I_{rk}, \quad  
\mathfrak{q} \mathcal{R} = \bigoplus_{k \geq 0} \mathfrak{q} I_{rk}. 
\) 
Now set 
\[ 
A = (\mathcal{R}' :_{\mathcal{R}} \mathfrak{p} \mathcal{R}) \quad \text{and} \quad  
B = \mathcal{R}' :_{\mathcal{R}} (\mathfrak{p} + \mathfrak{q}^{\infty}) \mathcal{R}. 
\] 
Then \(H = A / B\) is a finitely generated bigraded \(\mathcal{R}\)-module. Moreover,
\[ 
H_{(*,k)} = \frac{(I_{r(k+N+s+1)+i} : \mathfrak{p}) \cap I_{rk}}{(I_{r(k+N+s+1)+i} : \mathfrak{p} + \mathfrak{q}^{\infty}) \cap I_{rk}}. 
\] 
Now, by Proposition \ref{Proposition 2.8.},
\[
I_{r(k+N+s+1)+i}=I_{rN+i}I_r^{k+s+1}.
\]
Hence,
\[
(I_{r(k+N+s+1)+i}:\mathfrak{p})=(I_{rN+i}I_r^{k+s+1}:\mathfrak{p})\subseteq (I_r^{k+s+1}:\mathfrak{p})\subseteq (I_r^{k+s+1}:a).
\]
Since $a\in\mathfrak{p}$, the second containment follows. Again, using Lemma \ref{ArtinReeslemma} with $k+1$ in place of $k$, we have
\[
(I_r^{k+s+1}:a)\subseteq I_r^{k+1}\subseteq I_r^k=I_{rk}.
\]
Therefore,
\[ 
H_{(*,k)} = \frac{(I_{r(k+N+s+1)+i} : \mathfrak{p})}{(I_{r(k+N+s+1)+i} : \mathfrak{p} + \mathfrak{q}^{\infty})}. 
\] 
for all \(k \gg 0\). By applying \cite[Lemma 1.2]{Conca24}, we get 
\[ 
\mathrm{v}_{\mathfrak{p}}(I_{rk+i}) = \alpha(H_{(*,k-N-s-1)}). 
\] 
A straightforward modification of the arguments used in the proofs of \cite[Lemma 3.3]{BrunsConcaVarbaro} and \cite[Theorem 8.3.4]{BrunsConcaRaicuVarbaro2022} shows that \(\mathrm{v}_{\mathfrak{p}}(I_{rk+i})\) is eventually linear in \(k\). 

\end{proof}

\begin{corollary}\label{v(Ik)<reg(Ik)forhomogeneousideals}Assume, in addition, that $R$ is standard graded.
Let $\mathcal{I}=\{I_k\}_{k\geq 0}$ be a Noetherian graded family of homogeneous ideals in  $R$ such that, for all $k\geq 0$, we have
\(
I_{rk}=I_r^k.
\)
If
\(
\alpha(I_r)<\rho_R(I_r),
\)
then for large $k$,
\[
\mathrm{v}(I_k)<\operatorname{reg}(I_k).
\]
\end{corollary}

\begin{proof}
It follows from Theorem~\ref{Proposition 1.4}, Theorem \ref{reg/vnumberisquasilinear}, and Theorem \ref{vnumberisquasilinearforhomogeneous}.
\end{proof}

\begin{corollary}\label{v(Ik)<reg(Ik)formonomialideals}
    Let $\mathcal{I} = \{ I_k \}_{k \geq 0}$ be a Noetherian graded family of monomial ideals in $S$. If there exist at least two vertices of the Newton--Okounkov region $\Delta(\mathcal{I})$ with different coordinate sums, then for large $k$,
    \[\mathrm{v}(I_k) < \operatorname{reg}(I_k).\]
\end{corollary}

\begin{proof}

From  \cite[Theorem 5.1]{HaNguyenNguyen25}, $\lim\limits_{k\to\infty}\frac{\mathrm{reg}(I_k)}{k}=\delta(\Delta(\mathcal{I})) = \max\{ |\mathbf{v}| : \mathbf{v} \text{ is a vertex of } \Delta(\mathcal{I}) \}$ and from Theorem \ref{analogue} we have $\lim\limits_{k\to\infty}\frac{\mathrm{v}(I_k)}{k}=\lambda(\Delta(\mathcal{I})) = \min\{ |\mathbf{v}| : \mathbf{v} \text{ is a vertex of } \Delta(\mathcal{I}) \}$. This completes the proof.

\end{proof}

This provides a partial answer to \cite[Question 5.1]{FicarraSgroi24} of Ficarra and Sgroi.

\begin{definition}

 For a monomial \( x^{A} \in S \), write \( \mu(x^{A}) = \max\{ i : x_i \mid x^{A} \} \). A monomial ideal \( I \subseteq S \) is called stable if for every \( x^{A} \in I \) and every \( i < \mu(x^{A}) \) we have \( x_i \cdot \frac{x^{A}}{x_{\mu(x^{A})}} \in I \). It is called strongly stable if for every \( x^{A} \in I \) and every \( j > i \) with \( x_j \mid x^{A} \) we have \( x_i \cdot \frac{x^{A}}{x_j} \in I \).

\end{definition}

\begin{proposition}\label{thm:stable-quotient}
Let  $I \subseteq S$ be a stable monomial ideal, and let \( G(I) \) denote its unique minimal monomial generating set. Let $t= \max\{ \mu(u) \mid u \in G(I)\}$. Then\\
$(i)$ $P=(x_1, \ldots, x_t)\in \Ass(I)$ and\\
$(ii)$ $\mathrm{v}_P(I) = \min\{\deg(f) \mid f \in G(I) \text{ and } \mu(f)=t\} - 1$.

\end{proposition}

\begin{proof}
$(i)$ Let $ f = m/x_t $ where $m\in G(I)$ such that $\mu(m)=t$. 
The inclusion \( P \subseteq (I : f) \) follows directly from the definition of stable monomial ideal.  

For the reverse inclusion, take a monomial \( g \) such that \( g \in (I : f) \).  
Assume for contradiction that \( g \notin P \), then \( \operatorname{supp}(g) \subseteq \{x_{t+1},\dots,x_n\} \), where \(\operatorname{supp}(g)\) denotes the set of variables dividing \(g\). 
Since \( I \) is a monomial ideal, there exists \( u \in G(I) \) such that \( u \mid gf \).  
Write \( u = u_1u_2 \) where \( u_1 = \gcd(u,g) \) and \( u_2 = \gcd(u,f) \).  
Because \( \operatorname{supp}(g) \cap \operatorname{supp}(f) = \emptyset \), the factors \( u_1 \) and \( u_2 \) are coprime, with \( u_1 \mid g \) and \( u_2 \mid f \).
If \( u_1 \neq 1 \), then let j = $\min\{ i \mid x_i \text{ divides } u_1 \}$.  
Since \( \operatorname{supp}(u_1) \subseteq \operatorname{supp}(g) \subseteq \{x_{t+1},\dots,x_n\} \), we have \( j > t \).  
Therefore \( \mu(u) \ge j > t \), contradicting the choice of \( t \) as the maximum of \( \mu(u) \) over all \( u \in G(I) \).
Hence \( u_1 = 1 \) and \( u = u_2 \mid f \).  
Since \( f = m/x_t \), we have \( u \mid m \).  
Both \( u \) and \( m \) are minimal generators of \( I \), so \( u = m \).  
But then \( m \mid f = m/x_t \), which is a contradiction.
Thus \( g \) $\in$ \( P \), which implies \( (I:f) \subseteq P \). This completes $(i)$.  \\
$(ii)$ From the proof of part $(i)$, we have seen that \(P = (I:f)\) where \(f = m/x_t\) for some \(m \in G(I)\) with \(\mu(m)=t\). Let \(m' \in G(I)\) be a generator of minimal degree among those satisfying \(\mu(m')=t\). Then it is clear that \(\mathrm{v}_P(I) \le \deg(m') - 1\). To prove the reverse inequality, let \(f''\) be a monomial of minimal degree such that \((I:f'') = P\) and \(\mathrm{v}_P(I) = \deg(f'')\). Since \(x_t \in P\), we have \(x_t f'' \in I\). Hence there exists \(m'' \in G(I)\) dividing \(x_t f''\) with \(x_t f'' = m'' g\) for some monomial \(g\). If $x_t\mid g$, then $f''=m''(g/x_t)\in I$, a contradiction. Set \(h = m''/x_t\). From the definition of stable monomial ideal we have for  any \(i \le t\), \(x_i h \in I\), so \(P \subseteq (I:h)\). Conversely, if \(u h \in I\), then \(u f'' = u h g \in I\), which gives \(u \in (I:f'') = P\). Thus \((I:h) \subseteq P\), and consequently \((I:h) = P\). By the minimality of \(\deg(f'')\), we must have \(g = 1\). Therefore \(f'' = h = m''/x_t\). Now, since \(m''\) is a minimal generator with \(\mu(m'') = t\), we have \(\deg(m'') \ge \deg(m')\). Hence
\(
\mathrm{v}_P(I) = \deg(f'') = \deg(m'') - 1 \ge \deg(m') - 1.
\)
This completes the proof.
\end{proof}

 The following examples illustrate Proposition~\ref{thm:stable-quotient}. They highlight that the choice of generator matters, not every colon ideal obtained from a generator yields an associated prime and also demonstrate the difference between the local \(\mathrm{v}\)-number \(\mathrm{v}_P(I)\) and the  \(\mathrm{v}\)-number \(\mathrm{v}(I)\).

\begin{example}
Let \(I = (x^2, xy^2, y^3, y^2z)  \subseteq K[x,y,z]\) be a stable monomial ideal with respect to the order \( x < y <z \)  . 
We have \(\operatorname{Ass}(I) = \{ (x,y), (x,y,z) \}\) and \(\mathrm{v}(I) = 2\). 
Observe that \((I:(y^2z/z)) = (I:y^2) = (x,y,z)\) and \((I:(xy^2/y)) = (I:xy) = (x,y)\), whereas \((I:(x^2/x)) = (x,y^2) \notin \operatorname{Ass}(I)\).
\end{example}

\begin{example}

Let \(I = (x^2, xy, y^3) \subseteq K[x,y]\) be a stable monomial ideal with respect to the order \( x < y  \)  . 
We have \(\operatorname{Ass}(I) = (x,y)\) and \(\mathrm{v}(I) = 1\). 
Observe that \((I:(xy/y)) = (I:x) = (x,y)\) and \((I:(y^3/y)) = (I:y^2 )= (x,y)\). 
By Proposition \ref{thm:stable-quotient}, we take \(m = xy\) such that \(\mathrm{v}(I) = \deg m - 1\). Moreover, for \(I = (x^3, x^2y, xy^2, y^3, x^2z^2)\), a stable monomial ideal in \(K[x,y,z]\) with respect to the order \( x < y <z \), 
we have \(\operatorname{Ass}(I) = \{(x,y), (x,y,z)\}\) and \((I:(x^2z^2/z)) = (I:x^2z) = (x,y,z)\). 
Here \(\mathrm{v}_m(I) = 3\) where \(m = (x,y,z)\) but \(\mathrm{v}(I) = 2\). 
Also, \((I:(y^3/y)) = (I:y^2) = (x,y)\) and \((I:(xy^2/y)) = (I:xy) = (x,y)\). 
Thus \(\mathrm{v}(I) = \deg y^2 = \deg xy = 2\). This shows that the local \(\mathrm{v}\)-number of \(I\) at \(\mathfrak{m}\) is not equal to the \(\mathrm{v}\)-number of \(I\).

\end{example}

\begin{corollary}\label{cor:v-bounds}
Let $I$ be a stable monomial ideal in $S$. Then
\[
\alpha(I) - 1 \leq \mathrm{v}(I) \leq d(I) - 1.
\]
\end{corollary}

\begin{proof}

By \cite[Lemma 3.1]{BiswasMandalSaha24} and Proposition \ref{thm:stable-quotient} we have \(\alpha(I) - 1 \leq \mathrm{v}(I)  \leq d(I) - 1.\)

\end{proof}

\begin{remark}

By \cite{{HerzogHibi11}}, we have the following relations for conditions on a monomial ideal
\[
\text{universal lexsegment} \implies \text{lexsegment} \implies \text{strongly stable} \implies \text{stable}.
\]

\end{remark}

\begin{corollary}\label{v(I)<reg(I)forstablemonomialideals}

Let $I$ be a stable monomial ideal (strongly stable, lexsegment, universal lexsegment, respectively). Then \[\mathrm{v}(I) < \operatorname{reg}(I).\]

\end{corollary}

\begin{proof}

By \cite[Corollary 7.2.3]{{HerzogHibi11}} we have \(\operatorname{reg}(I) = \max \{ \deg(m) : m \in G(I) \}\) and from Corollary \ref{cor:v-bounds} we have \(\mathrm{v}(I) <\operatorname{reg}(I)\).

\end{proof}

\begin{proposition}\label{productofstablemonomialsisstable}
Let \( I \) and \( J \) be stable monomial ideals in \( S \). Then the product \( IJ \) is also a stable monomial ideal.
\end{proposition}

\begin{proof}
Let \( I,J \subset S \) be stable monomial ideals. Take \( w = uv \in IJ \) with \( u \in I, v \in J \), and let \( \mu(w) = \max\{i \mid x_i \text{ divides } w\} \). Then \(
\mu(w) = \max\{\mu(u), \mu(v)\}
\). Fix \( i < \mu(w) \). If \( \mu(u) \ge \mu(v) \), then \( \mu(w) = \mu(u) \) and it implies \( x_i \cdot \frac{u}{x_{\mu(u)}} \in I \), so \( x_i \cdot \frac{w}{x_{\mu(w)}} = \bigl( x_i \cdot \frac{u}{x_{\mu(u)}} \bigr) v \in IJ \). If \( \mu(v) > \mu(u) \), then \( \mu(w) = \mu(v) \) and it implies  \( x_i \cdot \frac{v}{x_{\mu(v)}} \in J \), so \( x_i \cdot \frac{w}{x_{\mu(w)}} = u\bigl( x_i \cdot \frac{v}{x_{\mu(v)}} \bigr) \in IJ \). Thus \( x_i(w/x_{\mu(w)}) \in IJ \) for all \( i < \mu(w) \), proving \( IJ \) is stable.
\end{proof}

\noindent

If $I\subseteq S$ is a stable monomial ideal generated in degree $d$. Then $I$ has linear quotients and hence, by \cite[Lemma~4.1]{ConcaHerzog03}, a $d$-linear resolution over $S$. Moreover, by Proposition~\ref{productofstablemonomialsisstable}, every power $I^k$ is stable and generated in degree $kd$. Hence, Corollary~\ref{cor:v-bounds} yields $\mathrm{v}(I^k)=kd-1$ for every $k\geq1$. Thus, \cite[Conjecture~2.6]{Ficarra24} holds for equigenerated stable monomial ideals.

\begin{corollary}
If \(I\) is a stable monomial ideal in $S$, then \(\mathrm{v}(I^k) < \operatorname{reg}(I^k)\) for all \(k \ge 1\).
\end{corollary}

\begin{proof}
This follows directly from Corollary \ref{v(I)<reg(I)forstablemonomialideals} and Proposition \ref{productofstablemonomialsisstable}.

\end{proof}

Let $I$ be a homogeneous ideal in $S$. The quotient ring $S/I$ is standard graded, and it inherits the grading from $S$ by $S_i/I_i$ for every $i \in \mathbb{N}$. We have $\dim_K(S_i/I_i) < \infty$ for all $i \in \mathbb{N}$. The generating function $i \mapsto \dim_K(S_i/I_i)$ is called the Hilbert function of $S/I$. 

\noindent The Hilbert series of $S/I$ is
\[\operatorname{Hilb}_{S/I}(t) = \sum_{i \in \mathbb{N}} \dim_K(S_i/I_i) t^i = \frac{Q(t)}{(1-t)^d}\]
with $Q(t) \in \mathbb{Z}[t]$ and $Q(1) \neq 0$, then $e(S/I) = Q(1)$ where $e(S/I)$ is the multiplicity.
If $d = 0$, then $e(S/I) = Q(1) = \sum_{i \in \mathbb{N}} \dim_K(S_i/I_i) = \ell(S/I)$, the length of $S/I$.

\begin{proposition}\label{vNumbervsHilbertMultiplicity}
    Let $I$ be a zero-dimensional homogeneous ideal in $S$. Then $\mathrm{v}(I) < e(S/I)$.
\end{proposition}

\begin{proof}
    Since $I$ is a zero-dimensional homogeneous ideal in $S$, it is $\mathfrak{m}$-primary. Let $d = \mathrm{v}(I)$. Then by definition there exists a homogeneous element $f \in S_d$ such that $(I:f) = \mathfrak{m}$ and $f \notin I$. Hence $\overline{f} = f + I \in S/I$ is non-zero. We can write $\overline{f}$ as a linear combination of monomials of degree $d$. Since $\overline{f}$ is non-zero, at least one such monomial, say $x_{i_1}\cdots x_{i_d}$, is non-zero in $S/I$. For each $j = 1,\ldots,d-1$, the monomial $x_{i_1}\cdots x_{i_j}$ is also non-zero in $S/I$, which implies $\mathrm{v}(I) < 1+d \leq \dim_k(S/I) = \ell(S/I)$. Thus, $\mathrm{v}(I) < e(S/I)$.
\end{proof}

The above proposition fails when we do not assume $I$ to be zero-dimensional.

\begin{example}
    
Let $S = K[x,y]$ and $I = (x^3, x^2y, xy^2)$. The Hilbert series of $ S/I$ is
\[\operatorname{Hilb}_{S/I}(t) = \frac{1+t+t^2-2t^3}{1-t},\]
where $Q(t) = 1+t+t^2-2t^3$ and $Q(1) = 1$. Thus $e(S/I) = 1$, while $\mathrm{v}(I) = 2$, so $\mathrm{v}(I) > e(S/I)$.\\

\end{example}

\section*{Acknowledgements}
\noindent 
The authors would like to thank Prof. T\`{a}i Huy H\`{a}  and  Dr. Thái Thành  Nguy\~{\^e}n
for their helpful discussions and valuable suggestions.


\begin{thebibliography}{50}



  




\bibitem{VanmathiSarkar25} {Vanmathi. A and P. Sarkar},\textit{ v-Numbers of symbolic power filtrations,} 
Collect. Math. (2025). https://doi.org/10.1007/s13348-025-00476-w



\bibitem{AmbhoreSahaSengupta24} {S. B. Ambhore, K. Saha, I. Sengupta,} \textit {The v-number of binomial edge ideals.} Acta Mathematica
Vietnamica (2024). https://doi.org/10.1007/s40306-024-00540-w




\bibitem{Barvinok02}{A. Barvinok,}\textit{ A course in convexity} volume 54 of Graduate Studies in
Mathematics. American Mathematical Soc., 2002. ISBN 978-0-8218-
2968-4. doi: 10.1090/gsm/054



 \bibitem{BiswasMandal24}{P. Biswas and M. Mandal}, \textit{ A study of v-number for some monomial ideals}, Collect. Math. 2024, https://doi.org/10.1007/s13348-024-00451-x.




\bibitem{BiswasMandalSaha24}{P. Biswas, M. Mandal and K. Saha}, \textit{Asymptotic behaviour and stability index of v-numbers of graded ideals,} Vietnam Journal of Mathematics, 2026.


\bibitem{BrunsConcaVarbaro} W. Bruns, A. Conca, and M. Varbaro, \textit{Castelnuovo-Mumford regularity and powers}, Commutative algebra, Springer, Cham, 2021, pp. 147–158, DOI 10.1007/978-3-030-89694-2.

\bibitem{BrunsConcaRaicuVarbaro2022} W. Bruns, A. Conca, C. Raicu, and M. Varbaro, \textit{Determinants, Gr\"obner bases and cohomology}, Springer Monographs in Mathematics, Springer, Cham, 2022, DOI 10.1007/978-3-031-05480-8. MR4627943.



\bibitem{CDF+}{J. Camarneiro, B. Drabkin, D. Fragoso, W. Frendreiss, D. Hoffman, A. Seceleanu, T.
Tang, and S. Yang}
,\textit{ Convex bodies and asymptotic invariants for powers of monomial ideals}, J. Pure Appl. Algebra 226 (2022), no. 10, Paper No. 107089, 21 pp, DOI
10.1016/j.jpaa.2022.107089. MR4397130







\bibitem{Conca24}{A. Conca}, \textit{ A note on the v-invariant}, Proc. Am. Math. Soc. 152 (2024), no. 6, 2349–2351.

\bibitem{ConcaHerzog03}
A. Conca and J. Herzog,
\textit{Castelnuovo--Mumford regularity of products of ideals},
Collect. Math. \textbf{54} (2003), no.~2, 137--152.

\bibitem{Cooper20} {S. M. Cooper, A. Seceleanu, S. O. Tohaneanu, M. Vaz Pinto and R. H. Villarreal}, \textit{ Generalized minimum
distance functions and algebraic invariants of Geramita ideals}, Adv. Appl. Math.  112 (2020), 101940,
34 pp.




\bibitem{Cutkosky13}{ S.D. Cutkosky}, \textit{Multiplicities associated to graded families of ideals}, Algebra Number Theory 7 (2013),
2059–2083.



 \bibitem{Cutkosky14}{ S.D. Cutkosky}, \textit {Asymptotic multiplicities of graded families of ideals and linear series}, Adv. Math. 264
(2014), 55–113.



 \bibitem{CutkoskyHerzogTrung99}{S.D. Cutkosky, J. Herzog, and N.V. Trung}, \textit{Asymptotic behaviour of the Castelnuovo–Mumford regularity}, Compos. Math. 118, 243–261 (1999).


\bibitem{Ficarra26}{A. Ficarra}, {\textit{The stable set of associated primes of a complementary edge ideal}}, (2026), preprint 
https://doi.org/10.48550/arXiv.2603.02358


\bibitem{Ficarra24}{A. Ficarra}, {\textit{Simon conjecture and the  $\mathrm{v}$-number of monomial ideals}}, Collect. Math.(2024), https://doi.org/10.1007/s13348-024-00441-z.


\bibitem{FicarraMarques25}{A. Ficarra,  P. Macias Marques}, {\textit{ The v-function of powers of sums of ideals}}, J Algebr Comb
62, 13 (2025). https://doi.org/10.1007/s10801-025-01429-z



\bibitem{FicarraSgroi23}{A. Ficarra, E. Sgroi}, \textit {Asymptotic behaviour of the v-number of homogeneous ideals}, J. Algebra \textbf{704} (2026), 273--297

\bibitem{FicarraSgroi24}{A. Ficarra, E. Sgroi}, \textit {Asymptotic behavior of integer programming and the v-function of a graded
filtration}. J. Algebra Appl. (2026). https://doi.org/10.1142/S0219498826502361

\bibitem{FiorindoGhosh24}{L. Fiorindo, D. Ghosh}, \textit{Asymptotic behaviour of Vasconcelos invariants for products and powers of graded ideals}, 2024, preprint https://doi.org/10.48550/arXiv.2401.17815

\bibitem{FiorindoGhosh}{L. Fiorindo, D. Ghosh}, \textit{ On the asymptotic behavior of the Vasconcelos invariant for graded modules} Nagoya Math. J. (2025), 15 pp., published online
https://doi.org/10.1017/nmj.2024.33.



\bibitem{GhoshNanduriPramanik26}{D. Ghosh, R. Nanduri, and S. Pramanik}. \textit{Asymptotic prime divisors and  Vasconcelos invariant}, (2026), arXiv preprint https://doi.org/10.48550/arXiv.2603.11961


\bibitem{GhoshPramanik25} {D. Ghosh and S. Pramanik}. \textit {Asymptotic v-numbers of graded (co)homology
modules involving powers of an ideal}. J. Algebra, 671:61–74, 2025.




\bibitem{HaNguyenNguyen25}
 {H.T. H\`a, H. D. Nguyen, T. T. Nguy\~{\^e}n}, \textit {Asymptotic regularity of graded families of ideals},
(2025), preprint https://arxiv.org/abs/2501.07710










\bibitem{HaNguyen}{H.T. H\`a and T.T. Nguy\~{\^e}n}, \textit{ Newton–Okounkov body, Rees algebra, and analytic spread of graded families
of monomial ideals}, Trans. Amer. Math. Soc. Ser. B 11 (2024), 1065–1097.


\bibitem{HerzogHibi11}{J. Herzog and T. Hibi}, {\textit {Monomial Ideals. Springer-Verlag London Limited, 2011}}


 \bibitem{HerzogHibiTrung07}{J. Herzog, T. Hibi, and N.V. Trung}, {\textit {Symbolic powers of monomial ideals and vertex cover algebras}},
Adv. Math. 210 (2007), 304–322.

\bibitem{JaramilloVillarreal21}{D. Jaramillo and R. H. Villarreal}, {\textit {The v-number of edge ideals}}, J. Combin. Theory Ser.A 177 (2021),
105310.








\bibitem{KavehKhovanskii12} {K. Kaveh and A. Khovanskii}, \textit{ Newton–Okounkov bodies, semigroups of integral points, graded algebras
and intersection theory}, Ann. of Math. (2), 176 (2012), no. 2, 925–978.


\bibitem{KavehKhovanskii14} {K. Kaveh and A. Khovanskii},{\textit{ Convex bodies and multiplicities of ideals}}, Proc. Steklov Inst. Math. 286,
no. 1 (2014), 268–284.


 \bibitem{Kodiyalam00}{V. Kodiyalam}, \textit {Asymptotic behaviour of Castelnuovo–Mumford regularity}, Proc. Amer. Math. Soc. 128,
no. 2 (2000), 407–411.




\bibitem{KotalSaha24} {N. Kotal, K. Saha}, {\textit {On the v-number of Gorenstein ideals and Frobenius powers}}, Bull. Malays.
Math. Sci. Soc. 47 (2024), 167, 17p.




\bibitem{KumarNanduriSaha25}{M. Kumar, R. Nanduri, K. Saha}, {\textit{The slope of the v-function and the Waldschmidt constant}}, J. Pure
Appl. Algebra 229 (2025), no. 2, 107881



\bibitem{LazarsfeldMustata09}
R.~Lazarsfeld and M.~Musta{\c{t}}\u{a},
\newblock \textit { Convex bodies associated to linear series,}
\newblock { Ann. Sci. \'Ec. Norm. Sup\'er. (4)},
{ 42} (2009), no.~5, 783--835.


 

 

\bibitem{Okounkov96} {A. Okounkov}, \textit  {Brunn–Minkowski inequality for multiplicities}, Invent. Math. 125 (1996), no. 3, 405–411.

 
 
\bibitem{Okounkov03} {A. Okounkov},\textit{ Why would multiplicities be log-concave? }In The orbit method in geometry and physics
(Marseille, 2000), volume 213 of Progr. Math., 329–347. Birkh\"auser Boston, Boston, MA, 2003.



\bibitem{Ratliff76} {L. J. Ratliff}, \textit{On prime divisors of \(I^n\),n-large}.  Michigan Math. J. 23 (1976)
337–352.



\bibitem{Saha24}{K. Saha}, \textit{ The v-Number and Castelnuovo-Mumford Regularity of Cover Ideals of Graphs}. Int.
Math. Res. Not. IMRN, no. 11, 9010–9019 (2024)




\bibitem{SahaSengupta22}{K. Saha, and I. Sengupta}. {\textit {The v-number of monomial ideals.}} J. Algebraic Combin.56,no.3(2022): 903–27.
https://doi.org/10.1007/s10801-022-01137-y.




\bibitem{HunekeSwanson06}{ I. Swanson and C. Huneke},{\textit{ Integral Closure of Ideals, Rings, and Modules, London Mathematical Society Lecture Note Series, Vol. 336, Cambridge University Press, Cambridge,
2006}}.







\bibitem{TrungWang}{N.V. Trung and H-J. Wang}, \textit{ On the asymptotic linearity of Castelnuovo–Mumford regularity}, J. Pure
Appl. Algebra 201 (2005), 42–48.







\end{thebibliography}
\end{document}